\def\titlep{C$^{*}$-bialgebra defined by the direct sum of Cuntz algebras
(The revised)}
\documentclass[11pt]{article}
\usepackage{graphicx,ifthen}
\usepackage{amssymb}

\font\germ=eufm10 at12pt

\def\goth#1{\hbox{\germ#1}}

\newcommand{\qed}{\hbox{\rule[-2pt]{3pt}{6pt}}}
\newcommand{\qedh}{\hfill\qed \\}

\newcommand{\vep}{\varepsilon}

\newtheorem{Thm}{Theorem}[section]
\newtheorem{ex}[Thm]{Example}
\newtheorem{defi}[Thm]{Definition}
\newtheorem{lemm}[Thm]{Lemma}

\newtheorem{rem}[Thm]{Remark}

\def\cal#1{\mathcal #1}
\def\con{{\cal O}_{n}}

\def\edot{=1,\ldots,n}
\def\pr{{\it Proof.}\quad}
\def\nset#1{\{1,\ldots,n\}^{#1}}

\def\co#1{{\cal O}_{#1}}
\addtocounter{footnote}{1}
\def\sftt#1{
\setcounter{equation}{0}
\addtocounter{footnote}{1}
\section{#1}
}

\def\ssft#1{\subsection{#1}}

\def\autherp{Katsunori Kawamura}
\def\emailp{{\it E-mail address}: kawamura@kurims.kyoto-u.ac.jp.}
\def\addressp{College of Science and Engineering Ritsumeikan University,\\
1-1-1 Noji Higashi, Kusatsu, Shiga 525-8577, Japan
}

\def\ptimes{\otimes_{\varphi}}
\def\delp{\Delta_{\varphi}}

\def\N{{\bf N}}
\def\tilco#1{\tilde{{\cal O}}_{#1}}
\def\hdel{\hat{\Delta}}
\begin{document}
%
%
\setcounter{section}{0}
\setcounter{footnote}{0}
\setcounter{page}{1}
\pagestyle{plain}

\title{\titlep}
\author{\autherp\thanks{\emailp}
\\
\addressp}
\maketitle

%
%
\begin{abstract}
We show that a tensor product among
representation of certain C$^{*}$-algebras induces a bialgebra.
Let $\tilde{{\cal O}}_{*}$ be the smallest unitization  
of the direct sum of Cuntz algebras
\[{\cal O}_{*}\equiv 
{\bf C}\oplus {\cal O}_{2}\oplus {\cal O}_{3}\oplus{\cal O}_{4}\oplus \cdots.\]
We show that there exists a non-cocommutative comultiplication 
$\Delta$ and a counit $\varepsilon$ of $\tilde{{\cal O}}_{*}$.
From $\Delta,\vep$ and the standard algebraic structure, 
$\tilde{{\cal O}}_{*}$ is a C$^{*}$-bialgebra.
Furthermore  we show the following: 
(i) The antipode on $\tilde{{\cal O}}_{*}$ never exist.
(ii) There exists a unique  Haar state  on $\tilde{{\cal O}}_{*}$.
(iii) For a certain one-parameter bialgebra automorphism group of 
$\tilde{{\cal O}}_{*}$, 
a KMS state on $\tilde{{\cal O}}_{*}$ exists.
\end{abstract}




%
%
\sftt{Introduction}
\label{section:first}
We study the representation theory of C$^{*}$-algebras.
It is often that the study of representations 
gives a new structure of algebras.
We show that a tensor product among
representation of C$^{*}$-algebras induces a bialgebra.
We start with the motivation of this study.
%
%
%
\ssft{Motivation}
\label{subsection:firstone}

Let $\con$ be the Cuntz algebra for $2\leq n<\infty$.
In \cite{TS01}, we introduced a tensor product $\ptimes$ among 
representations of Cuntz algebras:
%
%
\begin{equation}
\label{eqn:tensor}
{\rm Rep}\con\times {\rm Rep}\co{m}\ni
(\pi_{1},\pi_{2})\mapsto \pi_{1}\ptimes \pi_{2}\in{\rm Rep}\co{nm}
\end{equation}
for each $n,m\geq 2$ where we write ${\rm Rep}\con$ the class of 
unital $*$-representations of $\con$.
The tensor product $\ptimes$ is associative and distributive with respect 
to the direct sum but not symmetric, that is,
\[
\left\{
\begin{array}{l}
\left.
\begin{array}{l}
(\pi_{1}\ptimes\pi_{2})\ptimes \pi_{3}
=\pi_{1}\ptimes(\pi_{2}\ptimes \pi_{3}),\\
\\
\pi_{1}\ptimes(\pi_{2}\oplus  \pi_{3})=
\pi_{1}\ptimes \pi_{2}\oplus \pi_{1}\ptimes \pi_{3}\\
\end{array}
\right\}
\quad \mbox{ for any }\pi_{1},\pi_{2},\pi_{3},
\\
\\
\mbox{
\begin{minipage}[t]{10cm}
there exist $\pi_{1},\pi_{2}$ such that 
$\pi_{1}\ptimes \pi_{2}$ and   $\pi_{2}\ptimes \pi_{1}$
are not unitarily equivalent.
\end{minipage}
}
\end{array}
\right.\]
The decomposition formulae with respect to this tensor product
of irreducible representations are computed.
From this result, we inferred that there may exist
a bialgebra structure \cite{Kassel} associated with such tensor product.
On the other hand,
if an algebra $A$ has a counit $\vep$ which is an algebra morphism,
then the kernel of $\vep$ is a proper two-sided ideal of $A$.
From this, any simple algebra over a field $k$
never be a bialgebra except $A=k$.
In particular, any Cuntz algebra never be
because $\con$ is simple for each $n$.
%
%
\ssft{C$^{*}$-bialgebra}
\label{subsection:firsttwo}
The treatment of a bialgebra structure for a C$^{*}$-algebra $A$
is different from a purely algebraic case in several points \cite{KV,MD}.
For example, 
the definition of comultiplication of a C$^{*}$-algebra
is taken as a $*$-preserving linear map from  $A$ to
a certain completion of the algebraic tensor product $A\odot A$
because $A\odot A$ is not a C$^{*}$-algebra in general.
In this paper,  we define C$^{*}$-bialgebra as follows.
%
%
\begin{defi}
\label{defi:cstar}
A triplet $(A,\Delta,\vep)$ is a C$^{*}$-bialgebra
if $A$ is a unital C$^{*}$-algebra,
$\Delta$ is a unital $*$-homomorphism from
$A$ to the minimal tensor product $A\otimes A$ of $A$,
and $\vep$ is a unital $*$-homomorphism from $A$ to ${\bf C}$
which satisfy that 
%
%
\begin{equation}
\label{eqn:bialgebratwo}
(\Delta\otimes id)\circ \Delta=(id\otimes\Delta)\circ \Delta,\quad
(\vep\otimes id)\circ \Delta\cong id \cong (id\otimes \vep)\circ
\Delta.
\end{equation}
\end{defi}
We mention about relations among this C$^{*}$-bialgebra,
purely algebraic bialgebra and compact quantum group in
$\S$ \ref{subsection:secondtwo}.

%
%
\ssft{Main theorem}
\label{subsection:firstthree}
We introduce a C$^{*}$-bialgebra in order to explain the tensor product
$\ptimes$ in (\ref{eqn:tensor}).
We write the $1$-dimensional C$^{*}$-algebra by $\co{1}$
for convenience in this paper.
%
%
\begin{defi}
Let $\co{*}$ be the C$^{*}$-algebra consisting
of all functions $x$ from ${\bf N}$ to $\cup_{n\in {\bf N}}\con$ such that
$x_{n}\in \con$ for all $n$ and such that
$\|x_{n}\|\to 0$ as $n\to\infty$
with pointwise sum, product and involution
where ${\bf N}=\{1,2,3,\ldots\}$.
We call $\co{*}$ 
the {\it direct sum} of $\{\con:n\in {\bf N}\}$.
\end{defi}

\noindent
Remark that $\co{*}$ is a C$^{*}$-algebra without unit.
Hence we consider the smallest unitization $\tilco{*}$ of $\co{*}$,
that is, $\tilco{*}$ is a unital C$^{*}$-algebra with
a closed two-sided ideal $\co{*}$ such that
$\tilco{*}/\co{*}\cong {\bf C}$.

We state our main theorem as follows.
%
%
\begin{Thm}
\label{Thm:mainone}
There exists a data $(\Delta,\vep)$ such that
$(\tilco{*},\Delta,\vep)$ is a 
non-commutative and non-cocommutative C$^{*}$-bialgebra.
Furthermore the following holds:
\begin{enumerate}
\item
There exists a dense $*$-subalgebra ${\cal A}_{0}$ 
of $\tilco{*}$ with the common unit such that
$\Delta({\cal A}_{0})$ is included in the algebraic
tensor product ${\cal A}_{0}\odot {\cal A}_{0}$ of ${\cal A}_{0}$
and $({\cal A}_{0},\Delta|_{{\cal A}_{0}},\vep|_{{\cal A}_{0}})$ is a bialgebra
in the algebraic sense \cite{Kassel}.
\item
The image $\Delta$ of $\co{*}$ is contained in $\co{*}\otimes \co{*}$.
\item
The antipode of any dense subbialgebra of $(\tilco{*},\Delta,\vep)$ 
never exists.
\item
There exists a unique Haar state on $\tilco{*}$.
\item
For a certain 
one-parameter bialgebra automorphism group of $(\tilco{*},\Delta,\vep)$,
there exists a KMS state on $\tilco{*}$ \cite{BR}.
\item
The following holds for $\ptimes$ in (\ref{eqn:tensor}):
\[(\pi_{1}\otimes \pi_{2})\circ \Delta|_{\co{nm}}=\pi_{1}\ptimes
\pi_{2}\quad((\pi_{1},\pi_{2})\in{\rm Rep}\con\times {\rm Rep}\co{m})\]
where $\co{nm},\con$ and $\co{m}$
are naturally identified with subalgebras of $\tilco{*}$ and  
$\pi_{1}\otimes \pi_{2}$ is naturally identified
with the representation of $\tilco{*}\otimes \tilco{*}$.
\end{enumerate}
\end{Thm}

The C$^{*}$-algebra $\tilco{*}$ is abstract C$^{*}$-algebra which is not
defined on a Hilbert space.
It is a rare example that for every $n\geq 2$, $\con$'s appear all at once.
The extension of the direct sum of canonical endomorphisms of Cuntz algebras
is also a bialgebra endomorphism of $\tilco{*}$ 
($\S$ \ref{subsection:fourthtwo}).

The bialgebra $(\tilco{*},\Delta,\vep)$ is 
not a deformation of known algebra or group algebra.
There is no standard comultiplication of $\tilco{*}$.
The origin of the comultiplication of $\tilco{*}$ is a tensor product of 
representations of Cuntz algebras.

In $\S$ \ref{section:second},
we review basic definitions about
C$^{*}$-bialgebra, Cuntz algebras, direct sum of C$^{*}$-algebras and
smallest unitization.
In $\S$ \ref{section:third},
Theorem \ref{Thm:mainone} is proved.
In $\S$ \ref{section:fourth}, we show examples of
bialgebra endomorphisms and automorphisms of $\tilco{*}$.
In $\S$ \ref{section:fifth}, we show examples of C$^{*}$-subbialgebras
of $\tilco{*}$.
It is shown that the smallest unitization of 
the direct sum of matrix algebras
\[{\bf C}\oplus M_{2}({\bf C})\oplus M_{3}({\bf C})\oplus
M_{4}({\bf C})\oplus \cdots\]
is a C$^{*}$-subbialgebra of $\tilco{*}$ 
in $\S$ \ref{subsection:fifthfour}.

%
%
\sftt{Basic definitions and axioms}
\label{section:second}
%
%
%
\ssft{Bialgebra, antipode and morphism}
\label{subsection:secondone}
We review bialgebra according to \cite{Kassel}.
In this subsection,
any tensor product means the algebraic tensor product.
Let $k$ be the ground field with unit $1$.
A {\it coalgebra} is a triplet $(C,\Delta,\vep)$ where
$C$ is a vector space and $\Delta:C\to C\otimes C$
and $\vep:C\to k$ are linear maps satisfying 
the following axioms:
\[(\Delta\otimes id)\circ \Delta=(id\otimes \Delta)\circ \Delta,\quad
(\vep\otimes id)\circ \Delta\cong id\cong (id\otimes \vep)\circ \Delta.\]
A {\it bialgebra} is a quintuple $(B,m,\eta,\Delta,\vep)$
where $(B,m,\eta)$ is a unital associative algebra and 
$(B,\Delta,\vep)$ is a counital coassociative coalgebra
such that both $\Delta$ and $\vep$
are unital algebra morphisms.
For two bialgebras
$(B,m_{B},\eta_{B},\Delta_{B},\vep_{B})$ and 
$(A,m_{A},\eta_{A},\Delta_{A},\vep_{A})$ over $k$,
$f$ is a {\it bialgebra morphism} from $A$ to $B$ if
$f$ is an algebra morphism from 
$(A,m_{A},\eta_{A})$ to $(B,m_{B},\eta_{B})$ and 
$f$ is a coalgebra morphism from 
$(A,\Delta_{A},\vep_{A})$ to $(B,\Delta_{B},\vep_{B})$.
An endomorphism $S$ of $B$ is called an {\it antipode}
of $(B,m,\eta,\Delta,\vep)$ if $S$ satisfies
$m\circ (id\otimes S)\circ \Delta=\eta\circ \vep=
m\circ (S\otimes id)\circ \Delta$.
If an antipode exists on $B$, then it is unique.
If $(B,m,\eta,\Delta,\vep)$ has the antipode $S$, then
$(B,m,\eta,\Delta,\vep,S)$ is called a {\it Hopf algebra}.

%
%
\begin{lemm}
\label{lem:antipode}
Let $(B,m,\eta,\Delta,\vep)$ be a bialgebra.
Define the linear map $W$ on $B\otimes B$ by
%
%
\begin{equation}
\label{eqn:kactakesaki}
W(a\otimes b)\equiv \Delta(a)(I\otimes b).
\end{equation}
for $a,b\in B$.
If $W$ is not invertible, then 
the antipode of $(B,m,\eta,\Delta,\vep)$ never exists.
\end{lemm}
%
%
\pr
According to Remark 1.8 in \cite{MNW}, we show that
if the antipode exists, then $W$ is invertible.
Assume that $S$ is the antipode of $(B,m,\eta,\Delta,\vep)$.
Define the linear map $V$ on $B\otimes B$ by
$V(a\otimes b)\equiv (id\otimes S)(\Delta(a))(I\otimes b)$ for $a,b\in B$.
For $a,b\in B$, we see that
\[VW(a\otimes b)=(id\otimes m)(id\otimes S\otimes id)(X_{a,b})\]
where $X_{a,b}\equiv \{((\Delta\otimes id)\circ \Delta)(a)\}
\cdot (I\otimes I\otimes b)$.
By the coassociativity,
$X_{a,b}= \{((id\otimes \Delta)\circ \Delta)(a)\}\cdot (I\otimes I\otimes b)$.
Since $S(I)=I$,
\[VW(a\otimes b)
=[\{id\otimes (m\circ (S\otimes id)\circ \Delta)\}\circ \Delta](a)
\cdot (I\otimes b).\]
Because $\eta\circ \vep = m\circ (S\otimes id)\circ \Delta $
and $\{(id\otimes \vep )\circ \Delta\}(a)=a\otimes 1$,
we obtain $VW(a\otimes b)=a\otimes b$.
Hence $V=W^{-1}$.
Therefore the statement holds.
\qedh

\noindent
For any bialgebra $(B,m,\eta,\Delta,\vep)$,
$W$ in (\ref{eqn:kactakesaki}) satisfies the Pentagon equation
\[W_{12}W_{13}W_{23}=W_{23}W_{12}\]
where we use the leg numbering notation \cite{BS}.
The map $W$ is not invertible in general.

For two bialgebras $A$ and $B$, 
a map $f$ from $A$ to $B$ is a {\it bialgebra morphism} if 
$f$ is a unital algebra morphism and 
$\Delta_{B}\circ f= (f\otimes f)\circ \Delta_{A}$ and $\vep_{B}\circ f=\vep_{A}$.
A map $f$ is a {\it bialgebra endomorphism} of $A$ if $f$ is a bialgebra morphism
from $A$ to $A$.
A map $f$ is a {\it bialgebra isomorphism} if $f$ is a bialgebra morphism
and bijective.
In addition, if $A=B$, then
$f$ is called a {\it bialgebra automorphism} of $B$.

%
%
\ssft{C$^{*}$-bialgebra and quantum group}
\label{subsection:secondtwo}
When one considers the C$^{*}$-algebra version of bialgebra,
it is unavoidable to review C$^{*}$-algebraic approaches for quantum group
because the notion of C$^{*}$-bialgebra
has been already treated in the context of quantum groups. 
According to \cite{KV}, the popular topic of quantum groups can be approached
from two essentially different directions. 
The first is algebraic in nature by \cite{Drinfeld,Jimbo}.
The second approach is analytic in nature
as the generalization of Pontryagin
duality for abelian locally compact groups. 

According to $\S$ 1 in \cite{EV},
the study of quantum group in theory of operator algebra
was begun from duality theory of locally compact group.
The theory is constructed as
an analogy of group operator algebra.
After studies of concrete examples by Drinfel'd \cite{Drinfeld}
and Woronowicz \cite{W1,W2}, many nontrivial examples are discovered.
These examples in theory of operator algebras 
are also deformation of group operator algebras
or deformation of function algebra of homogeneous space.

On the other hand, our motivation is originated from
representation theory of Cuntz algebra \cite{TS01}.
Our bialgebra was found in computation of tensor product of representations.

We review the definition of compact quantum group by Woronowicz \cite{KV,W1,W2}. 
Let $A\otimes A$ means the minimal C$^{*}$-tensor product on $A$ and 
let $M(A \otimes A)$ be the multiplier algebra of
$A\otimes A$. A $*$-homomorphism $\Delta$ from $A$ to $M(A\otimes A)$ 
is called {\it non-degenerate} if $\Delta(A)(A\otimes A)$
is dense in $A\otimes A$.
A non-degenerate $*$-homomorphism $\Delta$ from 
$A$ to  $M(A\otimes A)$ is called {\it coassociative} 
if $(id\otimes \Delta)\circ \Delta = (\Delta \otimes  id)\circ \Delta$. 
Then we call $\Delta$ a {\it comultiplication} of $A$.
A {\it locally compact quantum semigroup} is a pair 
$(A,\Delta)$ of a $C^{*}$-algebra
$A$ and a comultiplication $\Delta$ of $A$ \cite{MD}.
The pair $(A,\Delta)$ in Definition \ref{defi:cstar} is
a locally compact quantum semigroup.
Furthermore, we may be able to call it a compact quantum monoid
because there exists a unit and a counit for this $(A,\Delta)$.
%
%
%
\begin{defi}
\label{defi:cqgtwo}\cite{MD,VanDaele,W2}
For a C$^{*}$-algebra $A$ with unit $I$,
a pair $(A, \Delta)$ is called a compact quantum group if 
$\Delta$ is a comultiplication of $A$ such that
$\Delta(A)\subset A\otimes A$, and 
the spaces $\Delta(A)(A \otimes I)$ and $\Delta(A)(I \otimes A)$
are dense in $A \otimes A$
where $S_{1}S_{2}$ means the linear span of 
the set $\{xy:x\in S_{1},\,y\in S_{2}\}$
for two subsets $S_{1},S_{2}$ of $A\otimes A$.
\end{defi}

\noindent
The density condition in Definition \ref{defi:cqgtwo}
is called the {\it cancellation property}
or the {\it cancellation law} of $(A,\Delta)$ \cite{KV,MNW}.
In Remark \ref{rem:cancell}, it is shown that
the C$^{*}$-bialgebra $(\tilco{*},\Delta,\vep)$ in
Theorem \ref{Thm:mainone} does not satisfy the cancellation law.
%
%
\begin{defi}\cite{KV}
\label{defi:quantum}
Let $A$ be a unital C$^{*}$-algebra with a comultiplication $\Delta$.
A state $\omega$ on $A$ is a {\it Haar state} on the pair $(A,\Delta)$ 
if $\omega$ satisfies the following:
%
%
\begin{equation}
\label{eqn:haar}
(\omega\otimes id)\circ \Delta=\omega(\cdot)I
=(id\otimes \omega)\circ \Delta.
\end{equation}
\end{defi}

\noindent
The pivotal results about compact quantum group
are the existence of a unique Haar state,
and the existence of the antipode on a certain dense subbialgebra.
Hence Theorem \ref{Thm:mainone} shows both
differences and similarities between compact quantum group
and our example $(\tilco{*},\Delta,\vep)$.

%
%
\ssft{Cuntz algebra}
\label{subsection:secondthree}
For $n\geq 2$, 
let $\con$ be the {\it Cuntz algebra} \cite{C}, that is,
a C$^{*}$-algebra which is universally generated by
generators $s_{1},\ldots,s_{n}$ satisfying
$s_{i}^{*}s_{j}=\delta_{ij}I$ for $i,j\edot$ and 
$\sum_{i=1}^{n}s_{i}s_{i}^{*}=I$
where $I$ is the unit of $\con$.
Because $\con$ is simple, 
any homomorphism from $\con$ to a C$^{*}$-algebra is injective.
If $t_{1},\ldots,t_{n}$ are elements of a unital C$^{*}$-algebra
${\cal A}$ such that
$t_{1},\ldots,t_{n}$ satisfy the relations of canonical generators of $\con$,
then the correspondence $s_{i}\mapsto t_{i}$ for $i\edot$
is uniquely extended to a unital $*$-embedding
of $\con$ into ${\cal A}$ from the uniqueness of $\con$.
Therefore we simply call such a correspondence 
among generators by an embedding of $\con$ into ${\cal A}$.
For other studies about the bialgebra associated with the Cuntz algebra,
see \cite{CPZ,DR,Roberts}.

%
%
\ssft{Direct sum and tensor product of C$^{*}$-algebras}
\label{subsection:secondfour}
We write ${\rm Hom}(A,B)$ the set of all $*$-homomorphisms from
C$^{*}$-algebras $A$ to $B$.
Let $\{A_{n}:n\in {\bf N}\}$  be a family of C$^{*}$-algebras.
The set of functions $x$ from ${\bf N}$ into $\cup_{n}A_{n}$ such that
$x_{n}\in A_{n}$ for each $n$ in ${\bf N}$ and such that
$\|x_{n}\|\to 0$ as $n\to \infty$, is a C$^{*}$-algebra
with pointwise sum, product and involution.
We write this C$^{*}$-algebra by $\oplus \{A_{n}:n\in {\bf N}\}$ and
call it the {\it direct sum} of the $A_{n}$'s ($\S$ 1.2.4. in \cite{Ped}).
Remark that there exist different
definitions of direct sum of C$^{*}$-algebras \cite{Sakai}.
The algebraic direct sum $\oplus_{alg}\{A_{n}:n\in {\bf N}\}$ 
of the $A_{n}$'s is a dense $*$-subalgebra of $\oplus\{A_{n}:n\in {\bf N}\}$.

For $n,m\in {\bf N}$,
we write the minimal tensor product of $A_{m}$ and $A_{l}$ by
$A_{m}\otimes A_{l}$.
Define 
%
%
\begin{equation}
\label{eqn:number}
{\cal N}_{n}\equiv \{(m,l)\in {\bf N}:ml=n\}.
\end{equation}
The C$^{*}$-algebra
$C_{n}\equiv \oplus\{A_{m}\otimes A_{l}:(m,l)\in {\cal N}_{n}\}$
is naturally identified with a C$^{*}$-subalgebra of 
$A_{*}\otimes A_{*}$
where $A_{*}\equiv \oplus\{A_{n}:n\in {\bf N}\}$.
Then the following holds:
\[A_{*}\otimes A_{*}=\oplus\{C_{n}:n\in {\bf N}\}.\]

Let $\{B_{n}:n\in {\bf N}\}$ be a family of C$^{*}$-algebras.
Let $\{f_{n}:n\in {\bf N}\}$ be a family of $*$-homomorphisms
such that $f_{n}\in {\rm Hom}(A_{n},B_{n})$ for each $n$.
Define the map $f$ from 
$A_{*}\equiv 
\oplus \{A_{n}:n\in {\bf N}\}$ to $B_{*}\equiv \oplus\{B_{n}:n\in {\bf N}\}$ by 
\[(f(x))_{n}\equiv f_{n}(x_{n})\quad (n\in {\bf N},\,
x\in \oplus\{A_{m}:m\in {\bf N}\}).\]
Then $f\in {\rm Hom}(A_{*},B_{*})$.
We write $f$ by $\oplus\{f_{n}:n\in {\bf N}\}$.
If $A_{n}$ is a C$^{*}$-subalgebra of $B_{n}$ for each $n$,
then  $A_{*}$ is naturally
identified with a C$^{*}$-subalgebra of $B_{*}$.

%
%
\ssft{Smallest unitization of C$^{*}$-algebra and C$^{*}$-bialgebra}
\label{subsection:secondfive}
According to $\S$ 2.1 in \cite{Wegge},
we review the smallest unitization of a C$^{*}$-algebra.
For a C$^{*}$-algebra $A$ without unit,
we write $\{(a,x):a\in {\bf C},x\in A\}$ by $\tilde{A}$.
Then $\tilde{A}$ is ${\bf C}\oplus A$ as a vector space.
With respect to the norm
\[\|(a,x)\|\equiv {\rm sup}\{\|ay+xy\|:y\in A,\,\|y\|=1\},\]
$\tilde{A}$ is a Banach space.
By the operation
$(a,x)(b,y)=(ab,ay+bx+xy)$ for $(a,x),(b,y)\in \tilde{A}$,
$\tilde{A}$ is a C$^{*}$-algebra with unit 
$I_{\tilde{A}}=(1,0)$.
By the natural embedding $\iota$ of $A$ into $\tilde{A}$,
$A$ is a closed two-sided ideal of $\tilde{A}$ such that
$\tilde{A}/A\cong {\bf C}$.
Define the map $j_{A}$ from $\widetilde{A\otimes A}$ to
$\tilde{A}\otimes \tilde{A}$ by
\[j_{A}((a,z))\equiv aI_{\tilde{A}}\otimes I_{\tilde{A}}+(\iota\otimes \iota)(z)
\quad(z\in A\otimes A,a\in {\bf C}).\]
Then we can verify that $j_{A}\in {\rm Hom}(\widetilde{A\otimes A},
\tilde{A}\otimes \tilde{A})$ such that
$j_{A}(I_{\widetilde{A\otimes A}})=I_{\tilde{A}}\otimes I_{\tilde{A}}$.
If $A_{0}$ is a dense $*$-subalgebra of $A$,
then $\tilde{A}_{0}\equiv \{(a,x):a\in {\bf C},\,x\in A_{0}\}$ 
is also a dense $*$-subalgebra of $\tilde{A}$.
For a state $\omega$ on $A$,
the state $\tilde{\omega}$ on $\tilde{A}$ is defined by
$\tilde{\omega}((a,x))\equiv a+\omega(x)$
for $(a,x)\in \tilde{A}$.

Let $A$ and $B$ be C$^{*}$-algebras without unit.
If $\phi\in {\rm Hom}(A,B)$,
then there exists $\tilde{\phi}\in{\rm Hom}(\tilde{A},\tilde{B})$
such that $\tilde{\phi}|_{A}=\phi$ and 
$\tilde{\phi}(I_{\tilde{A}})=I_{\tilde{B}}$.
For $f\in {\rm Hom}(A,B\otimes B)$, define 
%
%
\begin{equation}
\label{eqn:hatmap}
\hat{f}\equiv j_{B}\circ \tilde{f}
\in {\rm Hom}(\tilde{A},\tilde{B}\otimes \tilde{B}).
\end{equation}
By the identification of $B\otimes B$ with
a C$^{*}$-subalgebra of $\tilde{B}\otimes \tilde{B}$,
we see that $\hat{f}|_{A}=f$ and
$\hat{f}(I_{\tilde{A}})=I_{\tilde{B}}\otimes I_{\tilde{B}}$.

For a C$^{*}$-algebra $A$ without unit,
if $\Delta$ is a $*$-comultiplication of $A$ and 
$\Delta\in {\rm Hom}(A,A\otimes A)$,
then we obtain 
$\hat{\Delta}\in {\rm Hom}(\tilde{A},\tilde{A}\otimes \tilde{A})$
according to (\ref{eqn:hatmap}).
%
%
\begin{lemm}
\label{lem:extensionunit}
Let $A$ be a C$^{*}$-algebra without unit and assume that
there exist
$\Delta\in {\rm Hom}(A,A\otimes A)$ and $\vep\in {\rm Hom}(A,{\bf C})$
such that $(A,\Delta,\vep)$ satisfies (\ref{eqn:bialgebratwo}).
Then $(\tilde{A},\hat{\Delta},\tilde{\vep})$ is a C$^{*}$-bialgebra
where 
\[\tilde{\vep}((a,x))\equiv a+\vep(x)\quad((a,x)\in \tilde{A}).\]
\end{lemm}
%
%
\pr
By definition,
$\hat{\Delta}((a,x))=aI\otimes I
+(\iota\otimes \iota)(\Delta(x))$ for $(a,x)\in\tilde{A}$
where we denote $I_{\tilde{A}}$ by $I$.
Because
\[(\hdel\otimes id)\circ (\iota\otimes \iota)=
(\iota\otimes \iota\otimes \iota)\circ (\Delta\otimes id),\quad
(id\otimes \hdel)\circ (\iota\otimes \iota)=
(\iota\otimes \iota\otimes \iota)\circ (id \otimes \Delta),\]
we obtain that
\[
\begin{array}{rl}
\{(\hdel\otimes id)\circ \hdel\}((a,x))
=&
aI\otimes I\otimes I+\{(\iota\otimes \iota\otimes \iota)\circ 
(\Delta\otimes id)\circ \Delta\}(x),
\\
\\
\{(id\otimes \hdel)\circ \hdel\}((a,x))
=&
aI\otimes I\otimes I+\{(\iota\otimes \iota\otimes \iota)\circ 
(id\otimes \Delta)\circ \Delta\}(x)
\end{array}
\]
where we simply write identity maps on both $A$ and $\tilde{A}$ by 
the same symbol $id$.
From these and the coassociativity of $\Delta$ on $A$,
we obtain 
$(\hat{\Delta}\otimes id)\circ \hat{\Delta}
=(id\otimes \hat{\Delta})\circ \hat{\Delta}$ on $\tilde{A}$.
Hence $\hat{\Delta}$ is a comultiplication of $\tilde{A}$.

By definition, $\tilde{\vep}(I)=1$ and $\tilde{\vep}\circ \iota=\vep$.
From this,
\[
\{(\tilde{\vep}\otimes id)\circ \hat{\Delta}\}((a,x))
=aI+\{\iota \circ (\vep\otimes id)\circ \Delta\}(x)=aI+(0,x)=(a,x).
\]
In the same way, we obtain
$\{(id\otimes \tilde{\vep})\circ \hat{\Delta}\}((a,x))=(a,x)$. 
In this way,
we see that $(\tilde{\vep}\otimes id)\circ \hat{\Delta}\cong id \cong 
(id\otimes \tilde{\vep})\circ \hat{\Delta}$.
In consequence,
$(\tilde{A},\hat{\Delta},\tilde{\vep})$ is a C$^{*}$-bialgebra.
\qedh

For two C$^{*}$-algebras $A$ and $B$ without unit,
if both $A$ and $B$ satisfy the assumption in 
Lemma \ref{lem:extensionunit}
and $f\in {\rm Hom}(A,B)$ satisfies that 
$\Delta_{B}\circ f=(f\otimes f)\circ \Delta_{A}$ and $\vep_{B}\circ f=\vep_{A}$,
then we can verify that
$\tilde{f}$ is a C$^{*}$-bialgebra morphism
from $(\tilde{A},\hat{\Delta}_{A},\tilde{\vep}_{A})$
to $(\tilde{B},\hat{\Delta}_{B},\tilde{\vep}_{B})$.

%
%
\sftt{Proof of Theorem \ref{Thm:mainone}}
\label{section:third}

\noindent
{\it Proof of Theorem \ref{Thm:mainone}.}
The algebraic direct sum of 
$\{\con\otimes \co{m}:1\leq n,m<\infty\}$ is dense in $\co{*}\otimes \co{*}$.
In order to define a coassociative comultiplication on $\co{*}$,
we introduce families of embeddings among $\{\con\}_{n\geq 1}$.
Let $I_{n}$ be the unit of $\con$ and 
let $s_{1}^{(n)},\ldots,s_{n}^{(n)}$ be canonical generators of $\con$
for $n\geq 1$ where $s_{1}^{(1)}\equiv I_{1}$.
We identify $\con$ as a C$^{*}$-subalgebra of $\co{*}$ for each $n\geq 1$.
Define the $*$-subalgebra $\co{n,0}$ of $\con$
generated by $s_{1}^{(n)},\ldots,s_{n}^{(n)}$.
Then the algebraic direct sum $\co{*,0}$ of $\{\co{n,0}:1\leq n<\infty\}$
is a dense $*$-subalgebra of $\co{*}$.
In $\co{*}$, remark that $s_{i}^{(n)}s_{j}^{(m)}
=s_{i}^{(n)}(s_{j}^{(m)})^{*}=(s_{i}^{(n)})^{*}s_{j}^{(m)}=0$ when $n\ne m$. 
For $n,m\geq 1$, define $\varphi_{n,m}\in {\rm Hom}(\co{nm},
\con\otimes \co{m})$ by
%
%
\begin{equation}
\label{eqn:embedding}
\varphi_{n,m}(s_{m(i-1)+j}^{(nm)})\equiv s_{i}^{(n)}\otimes s_{j}^{(m)}
\quad(i\edot,\,j=1,\ldots,m).
\end{equation}
Then we can verify that the following diagram is commutative
for each $n,m,l\in {\bf N}$:

\noindent
%
\def\one{
\put(50,250){$\co{nml}$}
\put(150,280){\vector(3,1){200}}
\put(170,350){$\varphi_{n,ml}$}
\put(150,240){\vector(3,-1){200}}
\put(170,160){$\varphi_{nm,l}$}
\put(430,350){$\co{n}\otimes \co{ml}$}
\put(430,150){$\co{nm}\otimes \co{l}$}
\put(650,350){\vector(3,-1){200}}
\put(750,350){$id_{n}\otimes \varphi_{m,l}$}
\put(650,170){\vector(3,1){200}}
\put(750,160){$\varphi_{n,m}\otimes id_{l}$}
\put(880,250){$\co{n}\otimes \co{m}\otimes \co{l}$}
}
\thicklines
\setlength{\unitlength}{.1mm}
\begin{picture}(1000,300)(80,120)
\put(100,0){\one}
\end{picture}

\noindent
in other words,
%
%
\begin{equation}
\label{eqn:weak}
(\varphi_{n,m}\otimes id_{l})\circ \varphi_{nm,l}
=(id_{n}\otimes \varphi_{m,l})\circ \varphi_{n,ml}\quad(n,m,l\in {\bf N})
\end{equation}
where $id_{x}$ is the identity map on $\co{x}$ for $x=n,l$.
Let ${\cal N}_{n}$ be as in (\ref{eqn:number}).
For $C_{n}\equiv \oplus\{\co{m}\otimes \co{l}:(m,l)\in {\cal N}_{n}\}$,
define $\Delta_{\varphi}^{(n)}\in {\rm Hom}(\con,C_{n})$ by
%
%
\begin{equation}
\label{eqn:partone}
\Delta_{\varphi}^{(n)}(x)
\equiv \sum_{(m,l)\in {\cal N}_{n}}\varphi_{m,l}(x)\quad(x\in \con).
\end{equation}
For example,
%
%
\begin{equation}
\label{eqn:six}
\Delta_{\varphi}^{(6)}(s_{2}^{(6)})
=I_{1}\otimes s^{(6)}_{2}
+s_{1}^{(2)}\otimes s_{2}^{(3)}
+s_{1}^{(3)}\otimes s_{2}^{(2)}
+s_{2}^{(6)}\otimes I_{1}.
\end{equation}
Let $C_{*,0}$ be the algebraic direct sum of $\{C_{n}:1\leq n<\infty\}$. 
Because $C_{*,0}$ is the algebraic direct sum
of $\{\co{m}\otimes \co{l}:1\leq m,l<\infty \}$,
$C_{*,0}$ is a dense $*$-subalgebra of $\co{*}\otimes \co{*}$.
Define the map $\delp$ 
from $\co{*,0}$ to $C_{*,0}$ by
%
%
\begin{equation}
\label{eqn:comultiplication}
\Delta_{\varphi}\equiv \oplus\{\delp^{(n)}:n\in {\bf N}\}.
\end{equation}
From (\ref{eqn:weak}), the following holds for $x\in \con$:
\[
\begin{array}{rl}
\{(\Delta_{\varphi}\otimes id)\circ \Delta_{\varphi}\}(x)
=&
\sum_{m,l,k\geq 1,\,mlk=n}(\varphi_{m,l}\otimes id_{k})(\varphi_{ml,k}(x))\\
=&
\sum_{m,l,k\geq 1,\,mlk=n}(id_{m}\otimes \varphi_{l,k})(\varphi_{m,lk}(x))\\
=&\{(id\otimes \Delta_{\varphi})\circ \Delta_{\varphi}\}(x).
\end{array}
\]
Hence
$(\Delta_{\varphi}\otimes id)\circ \Delta_{\varphi}
=(id\otimes \Delta_{\varphi})\circ \Delta_{\varphi}$ on $\co{*,0}$.
Therefore $\Delta_{\varphi}$ is a $*$-comultiplication of $\co{*,0}$.
We denote simply the extension of $\Delta_{\varphi}$ 
on the whole of $\co{*}$ by the same symbol $\delp$. 
In this way, we obtain a $*$-comultiplication $\Delta_{\varphi}$ of $\co{*}$.

Define $\vep_{0}\in {\rm Hom}(\co{*},{\bf C})$ by
%
%
\begin{equation}
\label{eqn:counit}
\vep_{0}(x)\equiv 0   \quad (x\in \con,\,n\geq 2),\quad
\vep_{0}(x)\equiv x\quad (x\in \co{1}={\bf C}).
\end{equation}
For $x\in \con$, we see that
\[
\{(\vep_{0}\otimes id)\circ \delp\}(x)
= (\vep_{0}\otimes id)(\sum_{(m,l)\in {\cal N}_{n}} \varphi_{m,l}(x))
= \varphi_{1,n}(x)
=I_{1}\otimes x.
\]
In the same way, we obtain
$\{(id\otimes \vep_{0})\circ \delp\}(x)=x\otimes I_{1}$.
Hence
$(\vep_{0}\otimes id)\circ \delp\cong id\cong(id\otimes \vep_{0})\circ \delp$.

Because $(\co{*},\delp,\vep_{0})$ satisfies
the assumption in Lemma \ref{lem:extensionunit},
$(\tilco{*},\Delta,\vep)$ is a bialgebra
for $\Delta\equiv \hat{\Delta}_{\varphi}$
and $\vep\equiv \tilde{\vep}_{0}$.
By (\ref{eqn:six}), $\Delta_{\varphi}$ is non-cocommutative.
Hence the statement holds.

\noindent
(i) 
For the $*$-algebra $\co{*,0}$ in the proof of (i),
$\tilco{*,0}$ satisfies
the condition of ${\cal A}_{0}$ in the statement. 

\noindent
(ii) 
By definition of $\delp$, the statement holds.

\noindent
(iii) 
Assume that ${\cal A}$ is a dense subbialgebra of $\tilco{*}$.
Then $I_{1}\in {\cal A}\cap \co{1}$
and ${\cal A}\cap \con\ne \{0\}$ for each $n\geq 2$.
For the map $W$ in (\ref{eqn:kactakesaki}),
$W(x\otimes I_{1})=0$ when $x\in{\cal A}\cap \con$ for each $n\geq 2$.
Because $x\otimes I_{1}\in ({\cal A}\odot {\cal A})\setminus\{0\}$
when $x\ne 0$,
$W$ is not invertible on ${\cal A}$.
By Lemma \ref{lem:antipode}, the statement holds.

\noindent
(iv)
Assume that $\omega$ is a Haar state on $(\tilco{*},\Delta)$.
We write $I$ the unit of $\tilco{*}$.
Then
\[\{(\omega\otimes id)\circ \Delta\}(I_{1})=
(\omega\otimes id)(I_{1}\otimes I_{1})
=\omega(I_{1})I_{1}.\]
By (\ref{eqn:haar}), $\omega(I_{1})I_{1}=\omega(I_{1})I$.
This implies that $\omega(I_{1})=0$.

Next, for a prime number $p$,
\[\{(\omega\otimes id)\circ \Delta\}(I_{p})=
(\omega\otimes id)(I_{1}\otimes I_{p}\oplus I_{p}\otimes I_{1})
=\omega(I_{p})I_{1}.\]
From this and (\ref{eqn:haar}),
$\omega(I_{p})=0$ for all prime number $p$.

Assume that $\omega(I_{p_{1}\cdots p_{l}})=0$
for any $l$ prime numbers $p_{1},\ldots,p_{l}$.
Then for $l+1$ prime numbers $p_{1},\ldots,p_{l+1}$,
\[\{(\omega\otimes id)\circ \Delta\}(I_{p_{1}\cdots p_{l+1}})=
\sum_{(a,b)\in {\cal N}_{p_{1}\cdots p_{l+1}}}\omega(I_{a})I_{b}.\]
In the R.H.S., 
$\omega(I_{a})=0$ except $a=p_{1}\cdots p_{l+1}$ 
because $a$ is a divisor of $p_{1}\cdots p_{l+1}$.
Hence
$\{(\omega\otimes id)\circ \Delta\}(I_{p_{1}\cdots p_{l+1}})
=\omega(I_{p_{1}\cdots p_{l+1}}) I_{1}$.
From (\ref{eqn:haar}), we obtain $\omega(I_{p_{1}\cdots p_{l+1}})=0$.
Hence we see that $\omega(I_{n})=0$ for every $n\in {\bf N}$ by induction.
From this, $\omega|_{\co{*}}=0$.
On the other hand,
\[\{(\omega\otimes id)\circ \Delta\}(I)
=\omega(I)I=\{(id\otimes \omega)\circ \Delta\}(I).\]
Therefore $\omega((a,x))=a$ for each $(a,x)\in\tilco{*}$.
Hence the statement holds. 

\noindent
(v) 
For $n\in {\bf N}$, 
define the one-parameter automorphism group $\kappa^{(n)}$ of $\con$ by
%
%
\begin{equation}
\label{eqn:kappathree}
\kappa^{(n)}_{t}(s_{i}^{(n)})\equiv n^{-\sqrt{-1}t}s_{i}^{(n)}
\quad (i\edot,\,t\in {\bf R}).
\end{equation}
Define $\kappa_{t}^{(*)}\equiv\oplus \{\kappa_{t}^{(n)}:n\in {\bf N}\}$.
Then we can verify that $\tilde{\kappa}^{(*)}$ is a one-parameter bialgebra
automorphism group of $\tilco{*}$.
From Example 5.3.27 in \cite{BR},
there exists a unique $\kappa^{(n)}$-KMS state $\omega^{(n)}$ on $\con$.
For a sequence $b=(b_{n})_{n\geq 1}$ of non-negative
real numbers such that $\sum_{n\geq 1}b_{n}=1$,
define the state $\omega^{(*)}_{b}$ on $\co{*}$ by
\[\omega^{(*)}_{b}=\sum_{n\geq 1}b_{n}\omega^{(n)}\]
where $\omega^{(n)}$ is naturally identified with a state on $\co{*}$.
Then $\tilde{\omega}^{(*)}_{b}$ is a $\tilde{\kappa}^{(*)}$-KMS state
on $\tilco{*}$ for each $b=(b_{n})_{n\geq 1}$.

\noindent
(vi) 
For a representation $\pi_{i}$ of $\co{n_{i}}$ for $i=1,2$,
we can naturally identify $\pi_{i}$ as a representation of $\co{*}$.
Then 
\[\pi_{1}\ptimes \pi_{2}=(\pi_{1}\otimes \pi_{2})\circ \delp\]
defines a new representation $\pi_{1}\ptimes \pi_{2}$ of $\co{*}$.
By definition, $(\pi_{1}\ptimes \pi_{2})(x)=0$
if $x\in \oplus\{\co{m}:m\in {\bf N},\,m\ne n_{1}n_{2}\}$.
Therefore $\pi_{1}\ptimes \pi_{2}
=(\pi_{1}\otimes \pi_{2})\circ \varphi_{n_{1},n_{2}}$ on $\co{n_{1}n_{2}}$.
In this way, our construction of tensor product of 
representations in $\S$ 1 of \cite{TS01} is reconstructed.
\qedh

%
%
\begin{rem}
\label{rem:cancell}
{\rm
Remark that the pair $(\tilco{*},\Delta)$ is not a compact quantum group 
in Definition \ref{defi:quantum}
because it does not satisfy the cancellation law.
We give the proof here.
By definition,
\[\Delta(\tilco{*})(I\otimes \tilco{*})
={\bf C}I\otimes I\oplus 
I\otimes \iota(\co{*})\oplus 
 (\iota\otimes \iota)(\delp(\co{*}))
\oplus  (\iota\otimes \iota)(\delp(\co{*}))(I\otimes \iota(\co{*})).\]
Hence $x\otimes I\not\in 
\Delta(\tilco{*})(I\otimes \tilco{*})$ for any $x\in \iota(\co{*})\setminus
\{0\}$.
Therefore $\Delta(\tilco{*})(I\otimes \tilco{*})$
is not dense in $\tilco{*}\otimes \tilco{*}$.
Hence the pair $(\tilco{*},\Delta)$ 
does not satisfy the cancellation law.
}
\end{rem}

%
%
\sftt{Symmetry of $\tilco{*}$}
\label{section:fourth}
Here we use notations  $(\tilco{*},\Delta,\vep)$,
$\{\varphi_{n,m}\}_{n,m\geq 1}$,
$\co{*}, \delp,\vep_{0}$ in $\S$ \ref{section:third}.
Assume that $\{f_{n}\}_{n\geq 1}$ is a family of 
unital $*$-endomorphisms such that $f_{n}\in {\rm End}\con$ for each $n$.
For $f_{*}\equiv \oplus\{f_{n}:n\in {\bf N}\}$,
consider the following condition:
%
%
\begin{equation}
\label{eqn:morphismtwo}
\varphi_{n,m}\circ f_{nm}=(f_{n}\otimes f_{m})\circ \varphi_{n,m}
\quad(n,m\in {\bf N}).
\end{equation}
If (\ref{eqn:morphismtwo}) holds,
then $f_{*}$ is a $*$-endomorphism of $\co{*}$
such that $\delp\circ f_{*}=(f_{*}\otimes f_{*})\circ \delp$
and $\vep_{0}\circ f_{*}=\vep_{0}$.
Hence $\tilde{f}_{*}$ is a $*$-bialgebra endomorphism 
of $(\tilco{*},\Delta,\vep)$.
By using this,
we show examples of $*$-bialgebra endomorphism of $(\tilco{*},\Delta,\vep)$.

%
%
\subsection{Bialgebra automorphism arising from a family of unitary matrices}
\label{subsection:fourthone}
For a unitary $g=(g_{ij})_{i,j=1}^{n}\in U(n)$,
define the automorphism $\alpha_{g}^{(n)}$ of $\con$  by
\[\alpha_{g}^{(n)}(s_{i}^{(n)})\equiv \sum_{j=1}^{n}g_{ji}s_{j}^{(n)}
\quad(i\edot).\]
Let $g=\{g^{(n)}\}_{n\geq 1}$ be a family of unitary matrices such that
$g^{(n)}\in U(n)$ and
%
%
\begin{equation}
\label{eqn:unitariescondition}
g^{(1)}=1,\quad
(g^{(nm)})_{m(a-1)+b,m(i-1)+j}=(g^{(n)})_{a,i}\cdot (g^{(m)})_{b,j}
\end{equation}
for each $a,i\in \nset{}$, $b,j\in\{1,\ldots,m\}$ and $n,m\geq 1$.
Then the following holds for
the family $\{\varphi_{n,m}\}_{n,m\geq 1}$:
\[(\alpha_{g^{(n)}}^{(n)}\otimes \alpha_{g^{(m)}}^{(m)})
\circ \varphi_{n,m}=\varphi_{n,m}\circ \alpha_{g^{(nm)}}^{(nm)}
\quad(n,m\geq 1).\]
From this,
if (\ref{eqn:unitariescondition}) is satisfied, then
$\tilde{\alpha}^{(*)}_{g}$ is a 
bialgebra automorphism of $\tilco{*}$
where $\alpha^{(*)}_{g}=\oplus\{\alpha^{(n)}_{g^{(n)}}:n\in {\bf N}\}$.
We can verify that the set ${\bf G}$ of all families $\{g^{(n)}\}_{n\geq 1}$ 
which satisfies (\ref{eqn:unitariescondition})
is a subgroup of the direct product group $U(1)\times U(2)\times\cdots$.
Therefore $\tilde{\alpha}^{(*)}$ is an action of ${\bf G}$ on
the bialgebra $(\tilco{*},\Delta,\vep)$.

Let $\{\sigma^{(n)}\}_{n\geq 1}$ be a family of permutations such that
$\sigma^{(n)}\in {\goth S}_{n}$ for each $n$.
Define $g^{(n)}\in U(n)$ by
$(g^{(n)})_{ij}=\delta_{\sigma^{(n)}(j),i}$ for $i,j\edot$.
If $\{\sigma^{(n)}\}_{n\geq 1}$ satisfies that
\[\sigma^{(nm)}(m(i-1)+j)=m(\sigma^{(n)}(i)-1)+\sigma^{(m)}(j)
\quad(i\edot,\,j=1,\ldots,m)\]
for each $n,m\geq 1$,
then we see that 
$\{g^{(n)}\}_{n\geq 1}$ satisfies (\ref{eqn:unitariescondition}).
For the subgroup ${\bf S}$ of ${\bf G}$ which
consists of such families of unitaries associated with permutations,
we obtain the action of ${\bf S}$ on the bialgebra $(\tilco{*},\Delta,\vep)$.
%
%
\begin{ex}
\label{ex:kms}
{\rm
For $n\geq 1$, define $\zeta^{(n)}\in {\rm Aut}\con$ by
\[\zeta^{(n)}(s_{i}^{(n)})\equiv s_{n-i+1}^{(n)}\quad(i\edot).\]
Define $\zeta^{(*)}\equiv \oplus\{\zeta^{(n)}:n\in {\bf N}\}\in{\rm Aut}\co{*}$.
Then $\tilde{\zeta}^{(*)}$ is a bialgebra automorphism of $\tilco{*}$
such that $(\tilde{\zeta}^{(*)})^{2}=id$.
}
\end{ex}

%
%
\subsection{Canonical endomorphism of $\co{*}$}
\label{subsection:fourthtwo}
For $n\geq 2$, let $\rho^{(n)}$ be the canonical endomorphism of $\con$, that is,
\[\rho^{(n)}(x)\equiv s_{1}^{(n)}x(s_{1}^{(n)})^{*}+\cdots+
s_{n}^{(n)}x(s_{n}^{(n)})^{*}\quad(x\in \con).\]
Define $\rho^{(1)}\equiv id_{\co{1}}$.
Then we can verify that
$\varphi_{n,m}\circ \rho^{(nm)}
=(\rho^{(n)}\otimes \rho^{(m)})\circ \varphi_{n,m}$ for each $n,m\geq 1$.
Hence $\tilde{\rho}^{(*)}$ 
is a bialgebra endomorphism of $\tilco{*}$
where $\rho^{(*)}=\oplus\{\rho^{(n)}:n\in {\bf N}\}$.

%
%
\sftt{C$^{*}$-subbialgebras of $\tilco{*}$}
\label{section:fifth}
For a C$^{*}$-bialgebra $(A,\Delta,\vep)$,
$B$ is a {\it C$^{*}$-subbialgebra} of $(A,\Delta,\vep)$
if $B$ is a C$^{*}$-subalgebra of $A$ with the common unit,
and $(B,\Delta|_{B},\vep|_{B})$ is a C$^{*}$-bialgebra.
We regard that ${\bf N}$ is an abelian monoid, that is,
an abelian semigroup with unit with respect to
the product of natural numbers.
Let ${\cal S}$ be the set of all submonoids of the monoid ${\bf N}$.
%
%
%
\subsection{C$^{*}$-subbialgebra associated with submonoid of ${\bf N}$}
\label{subsection:fifthone}
For $H\in {\cal S}$, define the closed two-sided ideal $\co{*}(H)$ of $\co{*}$ by
%
%
\begin{equation}
\label{eqn:oh}
\co{*}(H)\equiv \oplus\{\con:n\in H\}.
\end{equation}
Define $\Delta_{\varphi^{(H)}}
\in {\rm Hom}(\,\co{*}(H),\,\co{*}(H)\otimes \co{*}(H)\,)$ by
\[\Delta_{\varphi^{(H)}}\equiv \oplus\{\Delta_{\varphi^{(H)}}^{(n)}:n\in H\},
\quad
\Delta_{\varphi^{(H)}}^{(n)}\equiv \sum_{(m,l)\in {\cal N}_{n}(H)}\varphi_{m,l}\]
where ${\cal N}_{n}(H)\equiv \{(m,l)\in H^{2}:ml=n\}$.
Then $(\widetilde{\co{*}(H)},\widehat{\Delta_{\varphi^{(H)}}},
\vep|_{\co{*}(H)})$ is also a C$^{*}$-bialgebra.
We see that $\Delta_{\varphi^{(H)}}=  \Delta_{\varphi}|_{{\cal O}_{*}(H)}$
if and only if there exist prime numbers $p_{1},\ldots,p_{r}$
such that $H=\{p_{1}^{k_{1}}\cdots p_{r}^{k_{r}}:
k_{1},\ldots,k_{r}=0,1,2,\ldots\}$.
In this case, $\widetilde{\co{*}(H)}$ is a C$^{*}$-subbialgebra of $\tilco{*}$.
%
%
\subsection{$UHF$ subbialgebra of $\tilco{*}$}
\label{subsection:fifthtwo}
Let $\gamma^{(n)}$ be the $U(1)$-gauge action on $\con$ defined by
$\gamma^{(n)}(s_{i}^{(n)})\equiv z s_{i}^{(n)}$
for $z\in U(1)$, $i\edot$ and $n\geq 2$.
We define
\[UHF_{n}\equiv \{x\in\con:\mbox{ for all }z\in U(1),\,
\gamma_{z}^{(n)}(x)=x\}.\]
We define $\gamma_{z}^{(1)}\equiv id$ on $\co{1}$ for $z\in U(1)$
and $UHF_{1}\equiv \co{1}={\bf C}I_{1}$.
Then we can verify that $\varphi_{n,m}|_{UHF_{nm}}\in
{\rm Hom}(UHF_{nm},UHF_{n}\otimes UHF_{m})$ for each $n,m\geq 1$.
Define the C$^{*}$-subalgebra $UHF_{*}$ of $\co{*}$ by
%
%
\begin{equation}
\label{eqn:uhftwo}
UHF_{*}\equiv \oplus \{UHF_{n}:n\in {\bf N}\}.
\end{equation}
Then $\Delta_{\varphi}|_{UHF_{*}}$ belongs to
${\rm Hom}(UHF_{*},UHF_{*}\otimes UHF_{*})$.
Because $\Delta_{\varphi}|_{UHF_{*}}$ is a comultiplication  
and $\vep|_{UHF_{*}}$ is a counit of $UHF_{*}$,
$\widetilde{UHF}_{*}$ is a C$^{*}$-subbialgebra of $\tilco{*}$.
For each $H\in {\cal S}$, define
%
%
\begin{equation}
\label{eqn:uhfh}
UHF_{*}(H)\equiv UHF_{*}\cap {\cal O}_{*}(H).
\end{equation}
Then $\widetilde{UHF}_{*}(H)$ is a C$^{*}$-subbialgebra of $\tilco{*}(H)$.
For $\kappa^{(*)}$ in (\ref{eqn:kappathree}),
the fixed-point subbialgebra $\tilco{*}^{\tilde{\kappa}^{(*)}}$ 
of $\tilco{*}$ by $\tilde{\kappa}^{(*)}$ is $\widetilde{UHF}_{*}$.
%
%
\subsection{Commutative C$^{*}$-subbialgebra associated with Cantor sets}
\label{subsection:fifththree}
For $n,l\geq 1$, define the finite-dimensional commutative
C$^{*}$-subalgebra ${\cal C}_{n,l}$ of $\con$ by
\[{\cal C}_{n,l}\equiv{\rm Lin}\langle\{s_{J}^{(n)}(s_{J}^{(n)})^{*}
:J\in \nset{l}\}\rangle\]
where $s_{J}^{(n)}=s_{j_{1}}^{(n)}\cdots s_{j_{l}}^{(n)}$ 
for $J=(j_{1},\ldots,j_{l})$.
We see that ${\cal C}_{n,l}\cong {\bf C}^{n^{l}}$.
Define the commutative C$^{*}$-subalgebra ${\cal C}_{n}$ of $\con$ 
by the inductive limit of $\{{\cal C}_{n,l}\}_{l\geq 1}$
with respect to the natural inclusions among $\{{\cal C}_{n,l}\}_{l\geq 1}$.
Then 
\[{\cal C}_{n}=\overline{{\rm Lin}\langle\{s_{J}^{(n)}(s_{J}^{(n)})^{*}:
J\in \nset{+}\}\rangle}\]
where $\nset{+}=\bigcup_{l\geq 1}\nset{l}$. 
Let $X_{n}=\nset{\infty}$ be the compact Hausdorff space 
of the infinite direct product of the discrete topological space $\nset{}$ and 
$C(X_{n})$ is the C$^{*}$-algebra 
of all continuous complex-valued functions on $X_{n}$.
Then ${\cal C}_{n}\cong C(X_{n})$ for each $n\geq 1$.
Especially, ${\cal C}_{1}\cong {\bf C}$.
Define the  commutative C$^{*}$-subalgebra ${\cal C}_{*}$ of $\co{*}$ by
\[{\cal C}_{*}=\oplus\{{\cal C}_{n}:n\in {\bf N}\}.\]
%
%
\begin{lemm}
\label{lem:cantor}
\begin{enumerate}
\item
The restriction $\delp|_{{\cal C}_{*}}$ of $\delp$
on ${\cal C}_{*}$ is also a comultiplication of ${\cal C}_{*}$.
\item
Identify $C(X_{a})\otimes C(X_{b})$ with $C(X_{a}\times X_{b})$.
Then $\varphi_{a,b}$ is identified with the map from $C(X_{ab})$ to
$C(X_{a}\times X_{b})$ such that
\[\varphi_{a,b}(f)(J,K)=f( J*K)\]
for $J=(j_{i})_{i\geq 1}\in X_{a}$ and $K=(k_{i})_{i\geq 1}\in X_{b}$
where $J*K\in X_{ab}$ is defined by
$J*K=(b(j_{i}-1)+k_{i})_{i\geq 1}$.
\end{enumerate}
\end{lemm}
%
%
\pr
(i)
We see that
$\varphi_{n,m}(s_{m(i-1)+j}^{(nm)}(s_{m(i-1)+j}^{(nm)})^{*})
=s_{i}^{(n)}(s_{i}^{(n)})^{*}\otimes s_{j}^{(m)}(s_{j}^{(m)})^{*}$
for each $i\edot$ and $j=1,\ldots,m$.
From this,
$\varphi_{n,m}({\cal C}_{nm})\subset 
{\cal C}_{n}\otimes {\cal C}_{m}$ for each $n,m\geq 1$.
By definition of $\delp$, 
$\delp({\cal C}_{*})\subset {\cal C}_{*}\otimes {\cal C}_{*}$.
Hence the statement holds.

\noindent
(ii)
For any $f\in C(X_{ab})$,
there exists a sequence $\{f_{l}\}_{l\geq 1}$
such that $f_{l}\in {\cal C}_{ab,l}$ and $f_{l}\to f$
as $l\to\infty$.
If $f_{l}=\sum_{R\in\{1,\ldots,ab\}^{l}}c_{R}s_{R}^{(ab)}(s_{R}^{(ab)})^{*}$,
then
\[\varphi_{a,b}(f_{l})=
\sum_{S\in \{1,\ldots,a\}^{l}}\sum_{T\in \{1,\ldots,b\}^{l}}c_{R * T}
s_{R}^{(a)}(s_{R}^{(a)})^{*}\otimes s_{T}^{(b)}(s_{T}^{(b)})^{*}\]
where $R*T=(b(r_{i}-1)+t_{i})_{i=1}^{l}$
for $R=(r_{i})_{i=1}^{l}$ and $T=(t_{i})_{i=1}^{l}$.
By identifying ${\cal C}_{a,l}$ with $C(\{1\ldots,a\}^{l})\cong 
({\bf C}^{a})^{l}$, we obtain
\[\varphi_{a,b}(f_{l})(R,T)=c_{R * T}=f(R*T).\]
By taking the limit $l\to \infty$, we have the statement.
\qedh

\noindent
Because $\vep|_{{\cal C}_{*}}$ is also a counit of ${\cal C}_{*}$,
$\tilde{{\cal C}}_{*}$ is a commutative, 
non-cocommutative C$^{*}$-subbialgebra of $\tilco{*}$ by (\ref{eqn:six}).

%
%
\subsection{C$^{*}$-subbialgebra defined by the direct sum 
of finite-dimensional algebras}
\label{subsection:fifthfour}
Identify the matrix C$^{*}$-algebra $M_{n^{l}}({\bf C})$ with
the C$^{*}$-subalgebra
\[{\rm Lin}\langle\{s_{J}^{(n)}(s_{K}^{(n)})^{*}:J,K\in\nset{l}\}\rangle\]
of $\con$ for $n, l\geq 1$.
We obtain the following C$^{*}$-subalgebras of $\co{*}$ for each $l\geq 1$:
\[M_{*^{l}}({\bf C})\equiv \oplus \{M_{n^{l}}({\bf C}):n\in {\bf N}\},\quad
\quad {\bf C}^{*^{l}} \equiv \oplus\{{\bf C}^{n^{l}}:n\in {\bf N}\}\]
where the later is the diagonal part of the former at each component and
${\bf C}^{*^{l}}=M_{*^{l}}({\bf C})\cap {\cal C}_{*}$.
Then both $\widetilde{M_{*^{l}}({\bf C})}$ and $\tilde{{\bf C}}^{*^{l}}$ 
are C$^{*}$-subbialgebras of $\tilco{*}$.
In particular, the smallest
unitizations of the following are C$^{*}$-subbialgebras of $\tilco{*}$:
\[M_{*}({\bf C})= \oplus\{M_{n}({\bf C}):n\in {\bf N}\},\quad
{\bf C}^{*} = \oplus\{{\bf C}^{n}:n\in {\bf N}\}.\]
Let $\{E_{ij}^{(c)}\}_{i,j=1}^{c}$ be the matrix unit of $M_{c}({\bf C})$
such that $E_{ij}^{(c)}=s_{i}^{(c)}(s_{j}^{(c)})^{*}$ for $i,j=1,\ldots,c$.
From (\ref{eqn:six}),
\[\delp(E_{2,2}^{(6)})=
I_{1}\otimes E_{2,2}^{(6)}+E_{1,1}^{(2)}\otimes E_{2,2}^{(3)}
+E_{1,1}^{(3)}\otimes E_{2,2}^{(2)}+E_{2,2}^{(6)}\otimes I_{1}. \]
Therefore both $M_{*}({\bf C})$ and ${\bf C}^{*}$ are non-cocommutative.

Let ${\cal A}$ be a C$^{*}$-subalgebra of $\co{*}$
generated by $\{I_{n}:n\geq 1\}$.
Then 
\[{\cal A}\cong c_{0}\cong C_{0}({\bf N}),\quad
\tilde{{\cal A}}\cong C({\bf N}\cup\{\infty\})\]
where $c_{0}$ is the C$^{*}$-algebra consisting
of sequences $(x_{n})_{n\in {\bf N}}$ of complex numbers such that 
$|x_{n}|\to 0$ as $n\to\infty$,
$C_{0}({\bf N})$ is the C$^{*}$-algebra consisting
of complex continuous functions on ${\bf N}$ vanishing at infinity and  
${\bf N}\cup\{\infty\}$ is the one-point compactification of 
the locally compact Hausdorff space $\N$ 
with respect to the discrete topology (Remark 2.1.8 in \cite{Wegge}).
The C$^{*}$-bialgebra $\tilde{{\cal A}}$ 
is a cocommutative C$^{*}$-subbialgebra of $\tilco{*}$.
\\

\noindent
{\bf Acknowledgement:}

The author would like to express my sincere thanks to Takeshi Nozawa
and Stefaan Vaes.
Stefaan Vaes points out errors of the first version of this paper.


\end{document}